\documentclass{article}
\usepackage[utf8]{inputenc}
\usepackage{biblatex}
\usepackage[margin=1in]{geometry}
\usepackage{calrsfs}
\usepackage{float}
\usepackage{xcolor}
\usepackage{amsmath, amssymb, amsthm, dsfont, todonotes}
\usepackage{bbm}
\newcommand{\id}{\mathbbm{1}}
\usepackage{commath}
\usepackage{mathtools}
\usepackage{braket}
\usepackage{algorithm}
\usepackage{algpseudocode}
\usepackage{subfiles}
\usepackage[capitalise]{cleveref}
\crefname{Subsection}{Subsection}{Subsections}

\makeatletter
\newtheorem*{rep@theorem}{\rep@title}
\newcommand{\newreptheorem}[2]{%
\newenvironment{rep#1}[1]{%
 \def\rep@title{#2 \ref{##1}}%
 \begin{rep@theorem}}%
 {\end{rep@theorem}}}
\makeatother

\newcommand{\R}{\ensuremath\mathbb{R}}

\newcommand{\vertiii}[1]{{\left\vert\kern-0.25ex\left\vert\kern-0.25ex\left\vert #1\right\vert\kern-0.25ex\right\vert\kern-0.25ex\right\vert}}

\newtheorem{theorem}{Theorem}
\newreptheorem{theorem}{Theorem}
\newtheorem{lemma}[theorem]{Lemma}
\newreptheorem{lemma}{Lemma}
\newtheorem{proposition}[theorem]{Proposition}
\newreptheorem{proposition}{Proposition}
\newtheorem{example}[theorem]{Example}
\newtheorem{corollary}[theorem]{Corollary}
\newreptheorem{corollary}{Corollary}
\newtheorem{assumptions}[theorem]{Assumptions}
\newreptheorem{assumptions}{Assumptions}
\newtheorem{remark}[theorem]{Remark}

\theoremstyle{definition}
\newtheorem{definition}[theorem]{Definition}
\newreptheorem{definition}{Definition}

\usepackage{authblk}
\title{Signature asymptotics, empirical processes, and optimal transport}
\author[1,2]{Thomas Cass}
\author[1]{Remy Messadene}
\author[1]{William F. Turner}
\affil[1]{Department of Mathematics, Imperial College London}
\affil[2]{Alan Turing Institute, London}

\date{\today}
\begin{document}

\maketitle




\begin{abstract}
     \par Rough path theory \cite{friz2020course} provides one with the notion of signature, a graded family of tensors which characterise, up to a negligible equivalence class, and ordered stream of vector-valued data. In this article, we lay down the theoretical foundations for a connection between signature asymptotics, the theory of empirical processes, and Wasserstein distances, opening up the landscape and toolkit of the second and third in the study of the first. Our main contribution is to show that the Hilbert-Schmidt norm of the signature can be reinterpreted as a statement about the asymptotic behaviour of Wasserstein distances between two independent empirical measures of samples from the same underlying distribution. In the setting studied here, these measures are derived from samples from a probability distribution which is directly determined by geometrical properties of the underlying path. The general question of rates of convergence for these objects has been studied in depth in the recent monograph of Bobkov and Ledoux \cite{bobkov2014one}. To illustrate this new connection, we show how the above main result can be used to prove a more general version of the original asymptotic theorem of Hambly and Lyons \cite{hambly2010uniqueness}. We conclude by providing an explicit way to compute that limit in terms of a second-order differential equation.
\end{abstract}



\section{Introduction}
\subsection{Previous work}
\par The mathematical notion of a path captures the concept of a continuous time-ordered sequence of values. These objects and their generalisations,
occur widely throughout both pure and applied mathematics. For example, the
analysis of the sample paths of a stochastic process forms
a significant part of stochastic analysis, while time series analysis is an
established tool in modern statistics. Abstract paths are inherently
infinite-dimensional objects, and it is desirable to seek low-dimensional
summaries which capture some features of interest. A mathematically-principled approach to effecting this has gained prominence in recent years and led several new developments in time-series analysis \cite{gyurko2013extracting, lyons2014rough, levin2013learning, morrill2020generalised}, machine learning \cite{chevyrev2016characteristic, ni2020conditional}, deep learning \cite{kidger2019deep, kidger2020neural} and more recently in kernel methods \cite{cochrane2021sk, lemercier2021distribution, cass2020computing, kiraly2019kernels}. This approach 
involves using the (path) signature transform which,  in distinction to traditional methods based on sampling, is rooted in capturing the path by understanding its effects on any smooth non-linear controlled differential system.  To be more precise, if $\gamma$ is a path of finite $1$-variation defined on the closed interval $[a,b]\subset \mathbb{R}$ into $\mathbb{R}^d$. Then, given a smooth collection of vector fields $\left\{
V_{i}:i=1,..,d\right\} $ on $%
\mathbb{R}
^{e}$, we can in some circumstances write the response $\alpha :\left[ a,b%
\right] \rightarrow 
\mathbb{R}
^{e}$ of the controlled differential equation 
\[
d\alpha _{t}=\sum_{i=1}^{d}V_{i}\left( \alpha _{t}\right) d\gamma _{t}^{i},%
\text{ started at }\alpha _{a}
\]
in terms of a convergent series of iterated integrals of $\gamma$; that is%
\[
\alpha _{b}-\alpha _{a}=\sum_{k=1}^{\infty }\sum_{i_{1}...i_{k}=1}^{e}\left.
V_{i_{1}}V_{i_{3}}...V_{i_{k}}\text{Id}\right\vert _{\alpha
_{a}}\int_{a<t_{1}<t_{2}<...<t_{i_{k}}<b}d\gamma _{t_{1}}^{i_{1}}..d\gamma
_{t_{k}}^{i_{k}},
\]%
where Id denotes the identity function on $%
\mathbb{R}
^{e}$.

\par Using the above as motivation, we recall that the \textbf{signature} $S(\gamma)$ of $\gamma$ is defined as the collection of all iterated integrals 
\begin{align}
    S(\gamma):=\bigg(1, S(\gamma)^1_{[a,b]}, S(\gamma)^2_{[a,b]}, ...\bigg) \in T((\mathbb{R}^d)):=\prod_{k=0}^{\infty}(\mathbb{R}^d)^{\otimes k}
\end{align}
where
\begin{equation}
    S(\gamma)^k:=\{S(\gamma)^{k; \,I}\}_{I \in J}
\end{equation}
and $J:=\{(i_1,...,i_k) |\;1\leq i_1 , \dots , i_k \leq d\} \subset \mathbb{N}^k$ is a set of multi-index and where
\begin{equation}
    S(\gamma)^{k; \,i_1,...,i_k}:=\bigg( \int_{a\leq t_1 \leq t_2 \leq .. \leq t_k \leq b} d\gamma_{t_1}^{i_1} \otimes d\gamma_{t_2}^{i_2} \otimes... \otimes d\gamma_{t_k}^{i_k}   \bigg) \in (\mathbb{R}^d)^{\otimes k}
\end{equation}
A key theorem of Hambly and Lyons \cite{hambly2010uniqueness} states that the map $\gamma \mapsto S(\gamma)$ is one-to-one up to an equivalence relation on the space of paths which is called tree-like equivalence \cite{boedihardjo2016signature}. In this way, the signature offers a top-down summary of $\gamma$ allowing one a practical and efficient representation of the curve \cite{levin2013learning}. This approach has several pleasant theoretical and computational consequences. For example, the signature transform satisfies a universality property in that any continuous function $f:J\rightarrow \mathbb{R}$ from a compact subset $J$ of signature features can be arbitrarily well approximated by a linear functional \cite{levin2013learning}. This result inspired the development of several new methods and paradigms in time-series analysis \cite{gyurko2013extracting, lyons2014rough, levin2013learning, morrill2020generalised}, machine learning \cite{chevyrev2016characteristic, ni2020conditional} and more recently deep learning \cite{kidger2019deep, kidger2020neural}. A notable use of the signature transform has been its application in the popular field of kernel methods where the so-called \textit{signature kernel} \cite{kiraly2019kernels}, consisting of the inner product between two signatures, is introduced. In addition to being backed up by a rich theory \cite{lyons2007differential, friz2020course}, working with the inner product of signature features has proved itself to be a promising and effective approach to many tasks \cite{cochrane2021sk, cass2020computing} and has achieved state of the art performance for some of them \cite{lemercier2021distribution}. This growing interest in the use of the signature has also brought into focus methods for recovering properties of the underlying path from the signature.

\par Some terms in the signature are explicitly relatable to properties of the original path, e.g. the increment and the area can be recovered from the terms of order $1$ and order $2$ respectively. Recovering more granular information on the path demands a more sophisticated approach. A rich stream of recent work has tackled the explicit reconstruction of a path from its signature, e.g. based on a unicity result of the signature for Brownian motion sample paths \cite{le2013stratonovich}, one can consider a polygonal approximation to Brownian paths \cite{le2013stratonovich}, diffusions \cite{geng2016inversion}, a large class of Gaussian processes \cite{boedihardjo2015uniqueness} and even some deterministic paths \cite{geng2017reconstruction} based on the signature features only.  These approaches fundamentally exploit the full signature representation (and not a truncated version of it) which may not be available in some cases. In parallel, there have been other approaches to reconstruction. In \cite{lyons2017hyperbolic}, the hyperbolic development of the signature is exploited to obtain an inversion scheme for piecewise linear paths. On the other hand, \cite{lyons2014inverting} proposed a symmetrization procedure on the signature to which leads to a reconstruction algorithm in some cases. Both approaches have the advantage of being implementable.

\par Another branch of investigation has been the recovery of broad features of a path using the asymptotics of its signature, or functions of terms of its signature. The study of the latter has been an active area of research for the last 10 years \cite{hambly2010uniqueness, boedihardjo2019non, chang2018super, boedihardjo2019tail, boedihardjo2020path}. Recall that if $\gamma$ is absolutely continuous its \textbf{length} is defined to be
\begin{equation}
    L(t-a):=L(\gamma)_{[a;t]}:=\int_{a}^t |\gamma'_s|ds,
\end{equation}
and denoted by $l:=L(b-a)$. Such a curve admits a \textbf{unit-speed parametrisation} $\rho:[0,l]\rightarrow[a,b]$ defined as
\begin{equation}
    \rho(s)=L^{-1}(s)
\end{equation}
such that the path $\gamma \circ \rho$ is a unit speed curve. Hambly and Lyons \cite{hambly2010uniqueness} initially showed that the arc-length of a unit-speed path can be recovered from the asymptotics of the norm of terms in the signature under a broad-class of norms.  To be more concrete, they proved that if $\gamma:\left[  0,l\right]  \rightarrow V$ is a continuously
differentiable unit-speed curve (where $V$ is a finite dimensional Banach space), then it holds that
\begin{equation}
\lim_{n\rightarrow\infty}\frac{n!\norm{ S(\gamma)^{n}}}{l^{n}}=1,\label{strong}%
\end{equation}
for a class of norms $\norm{\cdot}$ on the tensor product $V^{\otimes n}$ which includes the
projective tensor norm, but excludes the Hilbert-Schmidt norm in the case
where $V$ is endowed with an inner-product. For the same class of norms, this
also implies the weaker statement
\begin{equation}
\lim_{n\rightarrow\infty}\left(  n!\norm{ S(\gamma)^{n}} \right)  ^{1/n}=l.\label{weak}%
\end{equation}
A natural question is whether other properties of the original path can also be recovered in a similar fashion.

\par When $V$ is an inner-product space and the Hilbert-Schmidt norm is considered,
\cite{hambly2010uniqueness} also show that \cref{weak} holds, but that statement \cref{strong}
fails to be true in general: indeed under the assumption that $\gamma$ is
three times continuously differentiable the following result \footnote{We refer to the right-hand side term $c(\gamma)$ as the \textit{Hambly-Lyons limit}.} is proved,
\begin{equation}
\lim_{n\rightarrow\infty}\frac{n!\norm{ S(\gamma)^{n}} }{l^{n}}=c\left(  \gamma\right)  :\mathbb{=E}\left[  \exp\left(
\int_{0}^{l}\left(  B_{0,s}\right)  ^{2}\left\langle \gamma_{s}'%
,\gamma_{s}'''\right\rangle ds\right)  \right]  ^{1/2},\label{HL id}%
\end{equation}
where $\left(  B_{0,s}\right)  _{s\in\lbrack0,l]}$ is a Brownian bridge
which returns to zero at time $l$. It is easily seen that $c\left(  \gamma\right)<1$ unless $\gamma$ is a straight line.
\par More recent articles have focused on proving statements similar to \cref{weak}
under the fewest possible assumptions on $\gamma$ and on the tensor norms. In the article \cite{chang2018super} for example, it was proved that if $\gamma$ is continuous
and of finite $1-$variation then, under any reasonable tensor norm, we have
that%
\begin{equation}
\lim_{n\rightarrow\infty}\left(  n!\norm{ S(\gamma)^{n}} \right)  ^{1/n}=\sup_{n\geq1}\left(  n!\norm{ S(\gamma)^{n}} \right)  ^{1/n}>0\label{inf id}%
\end{equation}
provided that the sequence $\left\{  \norm{S(\gamma)^{n}}
:n=1,2,...,\infty\right\}  $ does not contain an infinite subsequence of zeros. Boedhardjo and Geng \cite{boedihardjo2019non} strengthened this result by proving that the existence of such a subsequence is equivalent to the underlying path being tree-like, see also \cite{chang2018corrigendum}. Taken together, these articles prove that for the identity in \cref{inf id} holds true for a wide class of continuous bounded variation paths.
It is conjectured that for a tree-reduced path $\gamma$ the limit is exactly
the length of $\gamma$, see \cite{chang2018corrigendum}. This question remains open at the time of writing.

\subsection{Contributions}
\par In this article, we contribute to the effort of recovering the original path from its signature by laying down a novel route. We do so by explicitly relating Hilbert-Schmidt norm of projected signatures with $p$-Wasserstein distances between discrete probability measures, allowing the study of the former using tools of the latter. These measures are characterised in terms of $\gamma$ only through an integral equation, making the contribution of the geometrical properties of $\gamma$ (such as its curvature) in the limit of the norm explicit. To ease notation, from this point forwards, we consider unit speed paths parameterised on $[0.1]$ rather than $[0,l]$.

\par The core insight of this connection originated when realising that a theorem by del Barrio, Giné, and Utze (see \cref{thm:bgu} below) can be recasted to re-express the Hambly-Lyons limit $c(\gamma)$ as the limit of a $2$-Wasserstein distance between empirical measures (\cref{sec:hambly_lyons_limit_is_limit_wasserstein}). Formally, for a twice-continuously-differentiable unit speed path $\gamma:[0,1]\rightarrow \mathbb{R}^d$ that is regular enough (see exact conditions in \cref{prop:hambly_lyons_limit_as_limit_wasserstein}), we construct and show the existence of a measure $\mu$ on $\mathbb{R}$, prescribed in terms of the following integral equation for its cumulative distribution function $F$,
\begin{align}
    F(t)=\int_{0}^{t}|\gamma''_{F(s)}|^{-1} ds \geq ct, \quad \lim_{s\rightarrow 0^+}F(s)=0, \quad \lim_{s\rightarrow b^-}F(s)=1,
\end{align}
for some constants $c,b \in \mathbb{R}_+$. When coupled with the B-G-U Theorem, we show that
\begin{align}
    c(\gamma)=\lim_{n\rightarrow \infty}\mathbb{E}\left[\exp( -nW_2^2(\mu_n,\mu))\right],
\end{align}
i.e. the Hambly-Lyons limit $c(\gamma)$ can be seen as the limit of $2$-Wasserstein distances between $\mu$ and an empirical version $\mu_n:=\sum_{i=1}^n \delta_{F^{-1}(U_{i})}$, where $\delta$ denote the Dirac delta distribution and where $\{U_{i}\}_{i\in\{1,...,n\}}$ is a sample of $n$ independent uniform random variables on $[0,1]$. 

\par In \cref{sec:generalising_hambly_lyons_limit}, motivated by the above insight, we derive relationships between the signature inner product $\langle S(\gamma)^n, S(\sigma)^n \rangle$ and a series of $p$-Wasserstein distances. As a first application, we re-derive a generalised version of the Hambly-Lyons Limit Theorem (\cref{subsec:generalising_hambly_lyons}) through the lens of discrete optimal transport by exploiting asymptotic results of Wasserstein distances between empirical measures \cite{bobkov2014one}.

\begin{reptheorem}{thm:generalisation_hambly_lyons}[Generalised Hambly-Lyons Limit Theorem]
Let $\gamma:[0,1]\rightarrow \mathbb{R}^d$ be a twice-continuously- differentiable unit-speed path such that the map $s \mapsto \left\vert\gamma_s''\right\vert$ is non-vanishing, and differentiable with bounded derivative. Then,
\begin{align}
    \lim_{n\rightarrow \infty} n!\norm{S(\gamma)^n}=\mathbb{E}\left[\exp\left(-\int_0^1 (B_{0,s})^2\left\vert\gamma_s''\right\vert^2ds\right)\right]^{1/2}.
\end{align}
\end{reptheorem}

\par \cref{subsection:computing_hambly_lyons} concludes this article by presenting a way to practically compute that limit through the solving of a second order distributional differential equation.

\section*{Acknowledgements}
The authors would like to thank Sergey G. Bobkov and Michel Ledoux for insightful comments regarding possible relaxations of the assumptions of \cref{thm:bgu}.

\section{The Hambly-Lyons limit and Wasserstein distances}\label{sec:hambly_lyons_limit_is_limit_wasserstein}
\par This section outlines the proof of \cref{thm:generalisation_hambly_lyons}. We do so by first recalling a theorem by del Barrio, Giné, and Utze (\cref{thm:bgu}  below) and observing that under some conditions, the former allows the rewriting of the Hambly-Lyons limit $c(\gamma)$ as a limit in terms of the 2-Wasserstein distance between empirical measures. The rest of this section will investigate the assumptions needed on $\gamma$ for these conditions to be fulfilled.

\subsection{The B-G-U Theorem}
\par We now present the above-mentioned theorem by del Barrio, Giné, and Utze.
\begin{definition}[$J_p$ functional and $I$-function $I(t)$]\label{def:j_p_functional and i_function}
\par Let $X$ be a non-constant random variable with law $\mu$. Suppose that $\mu$ has a density $f$ w.r.t Lebesgue measure and let $F$ be the associated distribution function. The \textbf{$J_p$ functional} is defined as
\begin{equation}
    J_p(\mu):=\int_{-\infty}^\infty \frac{[F(x)(1-F(x))]^{p/2}}{f(x)^{p-1}}dx \label{def:j_2_functional}
\end{equation}
for $p\in\mathbb{N}$. Moreover, if $F$ admits an absolutely continuous inverse on $(0,1)$ (or equivalently, by virtue of Proposition A.17 in \cite{bobkov2014one}, if $\mu$ is supported on an interval, finite or not, and the absolutely continuous component of $\mu$ has on that interval an almost everywhere positive density), define the \textbf{I-function} $I(t)$ for almost all $t \in (0,1)$ as
\begin{equation}
    I(t):=f(F^{-1}(t))\overset{(A.18)}{=}\frac{1}{(F^{-1})'(t)}. \label{def:integral_quantity}
\end{equation}
where the last equality exploited Proposition A.18 in \cite{bobkov2014one}.
\end{definition}
For the rest of this section, we denote by $\mu_n$ the empirical measure defined as
\begin{equation}
    \mu_n:=\sum_{i=1}^n \delta_{X_i}
\end{equation}
where $X_i, \; i\in\{1,2,..,n\}$ is an i.i.d. sample of random variables sampled according to $\mu$. 
\begin{theorem}[B-G-U Theorem; \cite{del2005asymptotics}]\label{thm:bgu}
\par Let $\mu$ be a measure supported on $(a,b)\subset\mathbb{R}$ such that it admits a density $f$. Assume further that $f$ is positive and differentiable and satisfies
\begin{equation}
    \sup_{a<x<b}\frac{F(x)(1-F(x))}{f(x)^2}|f'(x)|
    <\infty. \label{eq:cond_proba_distrib_func_bgu}
\end{equation}
and that $J_2(\mu)< \infty$. Denote by $F$ its distribution function. Then,
\begin{equation}
    nW_2^2(\mu_n,\mu) \rightarrow \int_0^1 \frac{(B_{0,t})^2}{I(t)^2}dt \label{eq:bgu_limit}
\end{equation}
as $n \rightarrow \infty$ weakly in $\mathbb{R}$ where $B_{0,t}$ is a Brownian bridge starting at $0$ and vanishing at $t=1$.
\end{theorem}

\begin{remark}[Origin of assumption \cref{eq:cond_proba_distrib_func_bgu}]\label{rem:quantile_L2}
Condition \cref{eq:cond_proba_distrib_func_bgu} has been a recurrent assumption in several asymptotic results involving quantile processes and ultimately goes back to Csorgo and Revesz \cite{csorgo1978strong}. To better understand its connection with the B-G-U \cref{thm:bgu}, recall the following identity
\begin{align}
    W_2^2(\mu_n,\mu)=\int_{0}^1\left(F^{-1}_n(s) - F^{-1}(s)\right)^2ds=\frac{1}{n}\int_{0}^1 \xi_n(s)^2ds \label{eq:trajectories_2_wasserstein_convergence}
\end{align}
with $\xi_n(s):=\sqrt{n}\left(F^{-1}_n(s) - F^{-1}(s)\right)$ denoting the quantile process. It is shown in \cite{del2005asymptotics} that $\xi_n$ converges weakly in $L^2(0,1)$ to $B_{0,t}/I(t)$. For a slightly more general object, the so-called normed sample quantile process $\rho_n$, it can be shown that requiring condition \cref{eq:cond_proba_distrib_func_bgu} allows one to asymptotically control $\rho_n$. More details can be found in \cite{csorgHo1983quantile}.
\end{remark}
The following corollary will form a key component of our generalisation of the Hambly-Lyons result.
\begin{corollary}\label{cor:BGU_empirical}
    Let $\mu$ be a measure supported on $(a,b)\subset\mathbb{R}$ which satisfies the assumptions of \cref{thm:bgu}, then
    \begin{equation}
        nW_2^2(\mu_n^1,\mu_n^2)\to 2\int_0^1 \frac{(B_{0,t})^2}{I(t)^2}dt
    \end{equation}
    weakly in $\mathbb{R}$ as $n\to\infty$, where $\mu_n^1$ and $\mu_n^2$ are independent copies of empirical measures from $\mu$.
\end{corollary}
\begin{proof}
    Let $\xi_n^1$ and $\xi_n^1$ denote the quantile processes of $\mu_n^1$ and $\mu_n^2$ defined in \cref{rem:quantile_L2}. Since $L^2(0,1)$ is separable, \cite[Theorem 2.8]{billingsley} implies that $\xi_n^1-\xi_n^2$ converges weakly in $L^2(0,1)$ to $(B^1_{0,t}-B^2_{0,t})/I(t)$, for independent Brownian bridges $B^1$ and $B^2$. Since the difference of independent Brownian bridges is itself a Brownian bridge with twice the variance, we can apply the Continuous Mapping Theorem and \cref{rem:quantile_L2} to conclude that
    \begin{equation*}
        nW_2^2(\mu_n^1,\mu_n^2)=\int_0^1(\xi_n^1(t)-\xi_n^2(t))^2dt\to 2\int_0^1\frac{(B_{0,t})^2}{I(t)^2}dt
    \end{equation*}
    weakly in $\mathbb{R}$.
\end{proof}

\par The existence of a bridge between signature asymptotics and the theory of empirical processes is hinted at when one considers a special instance of the B-G-U \cref{thm:bgu}. Indeed, when applied to a regular enough class of measures $\mu$ satisfying $I(t)=|\gamma_t''|^{-1}$, the right-hand side of B-G-U exactly coincides with the Hambly-Lyons limit. The following remark formalises this observation.

\begin{corollary}[Recovering the Hambly-Lyons limit]\label{cor:coinicde_hamlby_lyons_bgu}
Let $\gamma:[0,1]\rightarrow \mathbb{R}^d$ be a twice-continuously-differentiable unit-speed path with non-vanishing second derivative. Consider a probability measure $\mu$ on $\mathbb{R}$ supported on $(a,b)$ for some $a,b \in \mathbb{R}$ and having density $f$. Assume that the following four conditions are satisfied:
\begin{enumerate}
\item[A.] The density function $f$ is positive and differentiable,
\item[B.] Its associated $I$-function satisfies $I(t)=|\gamma_t''|^{-1}$ almost everywhere,
\item[C.] The distribution function $F$ and density function $f$ satisfy
\begin{equation}
    \sup_{a<x<b}\frac{F(x)(1-F(x))}{f(x)^2}|f'(x)|
    <\infty,
\end{equation}
\item[D.] $J_2(\mu)< \infty$.
\end{enumerate}
Then the B-G-U \cref{thm:bgu} implies that the Hambly-Lyons limit in $c(\gamma)$ can be rewritten as the limit of a $2$-Wasserstein distance, i.e.
\begin{align}
    c(\gamma)=\lim_{n\rightarrow \infty}\mathbb{E}\left[\exp\left( -nW_2^2(\mu_n,\mu)\right)\right]^{1/2}.
\end{align}
\end{corollary}
\begin{proof}
\par The assumptions for the B-G-U \cref{thm:bgu} are satisfied and the limit
\begin{align}
    \lim_{n\rightarrow \infty}nW_2^2(\mu_n,\mu) = \int_0^1 B_{0,t}^2|\gamma''|^2 dt < \infty
\end{align}
exists under the stated assumptions. If one further assumes that the assumptions of the Hambly-Lyons limit \cref{HL id} holds (i.e. that $\gamma$ is of class $C^3$), then the limits coincide as $\gamma$ is unit-speed and satisfies $\langle \gamma',\gamma''' \rangle=-|\gamma''|^2$.
\end{proof}

\subsection{Existence and characterisation of admissible measures}

\par The rest of this section will focus on characterising the paths, $\gamma$, for which a measure $\mu$ exists satisfying conditions $A$, $B$, $C$ and $D$ of \cref{cor:coinicde_hamlby_lyons_bgu}. First, we recast the determination of such $\mu$ into the solving of an integral equation in terms of $\gamma$ only. Thereafter, we formulate assumptions on $\gamma$ ensuring the existence of a solution to this integral equation, whose associated measure $\mu$ satisfies the conditions $A$, $B$, $C$ and $D$ of \cref{cor:coinicde_hamlby_lyons_bgu}.

\par We start by rewriting condition $B$ as an explicit condition on the cumulative distribution function $F$ associated to the measure $\mu$.

\begin{remark}[Reformulating condition $B$ as an integral equation for $F$]\label{rem:condition_b_as_an_integral equation}
\par Let $\gamma:[0,1]\rightarrow \mathbb{R}^d$ be a twice-continuously-differentiable unit-speed path with non-vanishing second derivative. If $F$ admits an absolutely continuous inverse on $(0,1)$ (which is the case whenever condition $A$ holds; see Proposition A.17 from \cite{bobkov2014one}), then condition $B$ holds if and only if $F$ satisfies the following integral equation,
\begin{equation}
    F(t)=\int_{a}^{t}|\gamma''_{F(s)}|^{-1}ds, \label{eq:fundamental_theorem_calcullus_f_and_F}
\end{equation}
for all $t\in(a,b)$ with boundary conditions $\lim_{s\rightarrow a^+}F(s)=0$, $\lim_{s\rightarrow b^-}F(s)=1$. To see this, first suppose that conditions $A$ and $B$ hold, then
\begin{equation*}
     \int_{a}^t I(F(s))ds \overset{\eqref{def:integral_quantity}}{=} \int_{a}^tf(s)ds = F(t) \overset{B}{=} \int_{a}^{t}|\gamma''_{F(s)}|^{-1}ds.
\end{equation*}
Assume now that condition $A$ and \cref{eq:fundamental_theorem_calcullus_f_and_F} hold, then the Lebesgue Differentiation Theorem implies that the density function associated to $F$ can be written as
\begin{equation}
    f(s)=|\gamma''_{F(s)}|^{-1}, \quad \text{a.e.} \label{eq:relationship_density_and_gamma}.
\end{equation}
And so
\begin{equation*}
    I(s)=f(F^{-1}(s))=|\gamma''_{F(F^{-1}(s))}|^{-1}=|\gamma''_{s}|^{-1}, \quad \text{a.e.}.
\end{equation*}

\end{remark}

\begin{lemma}\label{lem:framing_lemma}
\par Let $\gamma:[0,1] \rightarrow \mathbb{R}^d$ be a twice-continuously-differentiable unit-speed path with non-vanishing second-derivative. Assume the existence of two constants $a<b \in \mathbb{R}$ and a twice-differentiable function $F:(a,b)\rightarrow [0,1]$ satisfying
\begin{align}
    F(t)=\int_{a}^{t}|\gamma''_{F(s)}|^{-1} ds, \quad \lim_{s\rightarrow a^+}F(s)=0,\ \lim_{s\rightarrow b^-}F(s)=1, \label{eq:integral_equation_for_F}
\end{align}
with bounded second derivative $F''$. Then,
\begin{itemize}
    \item[(i)] $F$ is a cumulative distribution function with support on $(a,b)$ and admits an absolutely continuous inverse $F^{-1}$ on $(0,1)$.
    \item[(ii)] $F$ fulfills condition $A$,$B$,$C$, and $D$.
\end{itemize}
\end{lemma}

\begin{proof}

\par That $F$ is a cumulative distribution function corresponding to some measure $\mu$ supported on $(a,b)$ follows immediately from \cref{eq:integral_equation_for_F} and differentiability of $F$. By virtue of proposition A.17 in \cite{bobkov2014one}, $F$ admits an absolutely continuous inverse on $(0,1)$ if and only if $\mu$ is supported on an interval, finite or not, and the absolutely continuous component of $\mu$ has on that interval has an a.e. positive density (with respect to Lebesgue measure). We show that the latter statement holds. Indeed, since $\gamma''$ is non-vanishing, the Lebesgue Differentiation Theorem implies that the density function of the absolutely continuous component of the probability measure $\mu$ associated to $F$ is positive almost everywhere and satisfies
\begin{align}
    f(s)=|\gamma''_{F(s)}|^{-1}, \quad \text{a.e.}.
\end{align}
implying the existence of an absolutely continuous inverse $F^{-1}$. This concludes (i). 
\par Regarding point (ii), as $F$ is twice-differentiable, its underlying measure $\mu$ is absolutely continuous and $f$ is its density which is positive and differentiable, implying the fulfillment of condition $A$. Condition $B$ follows from \cref{rem:condition_b_as_an_integral equation}. Finally, since $\left\vert\gamma''\right\vert$ is bounded from below (since it is continuous and non-vanishing) and $F''$ is bounded from above, it follows that conditions $C$ and $D$ are both satisfied.
\end{proof}
\par The above result states that a twice-differentiable solution $F$ prescribed by the integral equation \cref{eq:integral_equation_for_F}, with the property that $F''$ is bounded, is the cumulative distribution function of a measure $\mu$ required for the application of \cref{cor:coinicde_hamlby_lyons_bgu}. We now determine the conditions on $\gamma$ for which the existence of such an $F$ is guaranteed. Let $\gamma$ be as in \cref{lem:framing_lemma}. Then the following observations can be made.

\begin{itemize}

    \item[1.] Since $\gamma$ has bounded curvature, i.e. there exists a constant $c>0$ such that $|\gamma''_s|\leq \frac{1}{c}\in\mathbb{R}_+$ for all $s\in[0,1]$ (or equivalently said, if the map $s \mapsto |\gamma''_s|^{-1}$ is bounded by below by $c$), then
    \begin{align}
        F(t)\geq ct \implies F(b)=1
    \end{align}
    for a constant $b\in[a,a+c]$ as $F$ is monotone increasing. 

    \item[2.] Additionally, if the map $s \mapsto |\gamma''_s|^{-1}$ is Lipschitz continuous, the Picard-Lindelöf Theorem ensures the existence of a unique differentiable function $F:[a,b]\rightarrow [0,1]$ such that
    \begin{align}
        F'(t)=|\gamma''_{F(t)}|^{-1}>c, \quad F(a)=0. \label{eq:ODE_F}
    \end{align}

    \item[3.] To ensure that $F$ is twice differentiable, requiring Lipschitz continuity on $s \mapsto |\gamma_s''|^{-1}$ is not enough. Indeed, the latter assumption only implies the differentiability of $F'$ \textit{almost everywhere} on any open subset of the definition domain by virtue of the Rademacher's Theorem (which is not sufficient as it can lead to the breaking of condition $A$). However, if we assume the map $s \mapsto |\gamma''_s|$ to be differentiable, then $F$ is guaranteed to be twice-differentiable for all $t\in(a,b)$ by the quotient rule and non-vanishing property of $\gamma''$.
    
    \item[4.] As it is now assumed that $s \mapsto |\gamma_s''|^{-1}$ is differentiable, the fundamental theorem of calculus implies that the probability density function $f$ can be exactly written as
    \begin{align}
        f(s)=|\gamma_{F(s)}''|^{-1}\overset{\eqref{eq:ODE_F}}{\geq} c_+, \quad s \in (a,b).
    \end{align}
    Additionally, since $\gamma''$ is non vanishing, the curvature $|\gamma''_{s}|$ is bounded from below, i.e. $|\gamma''_s|\geq \frac{1}{c_-}\in\mathbb{R}_+$ for all $s\in[0,1]$ (or equivalently said, if the map $s \mapsto |\gamma''_s|^{-1}$ is bounded above by $c_-$), then the ODE \cref{eq:ODE_F} implies
    \begin{align}
        F'(t)\leq c_-.
    \end{align}
    Assuming further that the derivative of the map $s\mapsto |\gamma_{s}''|$ is bounded implies the boundedness of $F''$ by the chain and quotient rules.
\end{itemize}

\par These observations are collected and combined with the \cref{lem:framing_lemma} in the following lemma.

\begin{lemma}[Condition on $\gamma$ for existence of admissible $F$]\label{lem:existence_sol_framing}
\par Let $\gamma:[0,1]\rightarrow \mathbb{R}^d$ be a twice-continuously-differentiable unit-speed path such that the map $s \mapsto |\gamma_s''|$ is non-vanishing, and is differentiable with bounded derivative. Then there exists constants $a,b\in\mathbb{R}$ twice-differentiable function $F:(a,b) \rightarrow [0,1]$ such that the integral equation and boundary conditions \cref{eq:integral_equation_for_F} are satisfied and the conditions of \cref{cor:coinicde_hamlby_lyons_bgu} are also satisfied.
\end{lemma}
\begin{remark}[Boundedness of the derivative of $|\gamma_s''|$ as a condition on the curvature]
Observe that the boundedness condition on $\left(|\gamma_s''|\right)'$ (which we recall is there to ensure the fulfillment of the boundary conditions in equation \cref{eq:integral_equation_for_F}) can be reformulated as a condition on the total curvature $T(t):=\int_{0}^t |\gamma_s''|ds$ of $\gamma$ as follows. For a $c_2\in\mathbb{R}_+$,
\begin{align}
    c_2 & \geq \left(|\gamma_s''|^{-1}\right)'=-\frac{|\gamma_s''|'}{|\gamma_s''|^2}  =- \bigg[\log|\gamma_s''|\bigg]'\cdot|\gamma''_s|^{-1} \\
    \iff|\gamma_s''|&\geq -\frac{1}{c_2}\bigg[\log|\gamma_s''|\bigg]'
\end{align}
\par Taking the integral from $0$ to $t$ gives
\begin{align}
    \frac{1}{c_2}\log\left(\frac{|\gamma_0''|}{|\gamma_t''|}\right) \leq  T(t),\quad \forall t \in [0,1],
\end{align}
\end{remark}

\par Finally, the direct application of the existence \cref{lem:existence_sol_framing} followed by \cref{thm:bgu} yields the following result.

\begin{proposition}[Hambly-Lyons limit as the limit of a Wasserstein distance]\label{prop:hambly_lyons_limit_as_limit_wasserstein}
\par Let $\gamma$ be as in \cref{lem:existence_sol_framing} and $F$ the solution to the integral equation \cref{eq:integral_equation_for_F}. Let $\mu$ be the measure associated with $F$ and let $X_1,...,X_n$ be a sample of $n$ i.i.d. random variables drawn from $\mu$. Let $\mu_n$ be its associated empirical measure, i.e.
\begin{align}
    \mu_n:=\frac{1}{n}\sum_{i=1}^n \delta_{X_i}
\end{align}
Then, the Hambly-Lyons limit in $c(\gamma)$ is the limit of a $2$-Wasserstein distance
\begin{align}
    c(\gamma)=\lim_{n\rightarrow \infty}\mathbb{E}\left[\exp\left( -nW_2^2(\mu_n,\mu)\right)\right]^{1/2}.
\end{align}
\end{proposition}

\section{Signature projections and Wasserstein distances}\label{sec:generalising_hambly_lyons_limit}
\par Because of the known connection between the Hambly-Lyons limit and the limit of the Hilbert-Schmidt tensor norm of projected signatures \cite{hambly2010uniqueness}, the insights developed in the previous section naturally lead one to ask whether the Hilbert-Schmidt tensor norm of projected signatures can be related to Wasserstein distances. This section answers this question positively. By using this relationship, we are able characterise a class of curves larger than the $C^3$ one originally considered in \cite{hambly2010uniqueness} that satisfies
\begin{align}
    \lim_{n\rightarrow \infty} n!\norm{S(\gamma)^n}= c(\gamma).
\end{align}

We proceed as follows.

\begin{enumerate}
    \item First, in \cref{subsec:statistical_expression_sig_inner_product}, we prove a technical augmentation of a lemma in \cite{hambly2010uniqueness} and then derive a probabilistic representation for the inner product of two signature terms.
    \item Once this is done, \cref{subsec:lower_b_sig_in_terms_of_wasserstein} will exploit the characterisation of the Wasserstein distances between empirical measures to relate the quantities derived in the first step to these Wasserstein distances and hence derive lower and upper bounds on $\norm{S(\gamma)}$ in terms of the former.
    \item By leveraging the results of the previous section, we generalise the Hambly-Lyons limit \cref{HL id} in \cref{subsec:generalising_hambly_lyons} and present the proof of \cref{thm:generalisation_hambly_lyons}.
    \item Finally, we show a practical way to compute the limit in \cref{subsection:computing_hambly_lyons} and illustrate it in a simple case.
    
\end{enumerate}

\subsection{A probabilistic expression for signature inner products}\label{subsec:statistical_expression_sig_inner_product}
\par In this subsection, we generalise Lemma 3.9 in \cite{hambly2010uniqueness} to the inner product between signatures before presenting a probabilistic formula in terms of the angles between the derivatives of the two underlying curves.

\begin{definition}[Uniform order statistics sample]\label{def:uniform_order_statistics_samples}
    Let $\{U_i\}_{i \in \{1,2...,n\}}$ and $\{V_i\}_{i \in \{1,2...,n\}}$ be two independent collections of $n$ i.i.d. uniform random variables in $[0,1]$. Consider the relabeling $\{U_{(i)}\}_{i \in \{1,2...,n\}}$ where \begin{equation}
        0 \leq U_{(1)} \leq U_{(2)} \leq ... \leq U_{(n)}\leq 1.
    \end{equation}
    and similarly for $\{V_{(i)}\}_{i \in \{1,2...,n\}}$.
    In the rest of this article, we will denote by $\{U_{(i)}\}_{i \in \{1,2...,n\}}$ and  $\{V_{(i)}\}_{i \in \{1,2...,n\}}$ two independent collections of \textbf{i.i.d. uniform order statistics} on $[0,1]$.
\end{definition}

\begin{lemma}[Generalisation of Lemma 3.9 in \cite{hambly2010uniqueness}] \label{lem:sig_and_order_stat_augmented}
Let $\gamma,\sigma:[0,1] \rightarrow \mathbb{R}^n$ be two $C^1$ unit-speed curves and let $\{U_{(i)}\}_{i \in \{1,2...,n\}}$ and $\{V_{(i)}\}_{i \in \{1,2...,n\}}$ be two i.i.d. uniform order statistics collections (\cref{def:uniform_order_statistics_samples}). Then,
\begin{align}
\langle S(\gamma)^n, S(\sigma)^n \rangle = \frac{1}{(n!)^2} \mathbb{E}\bigg[\prod_{i=1}^{n} \langle \gamma'_{U_{(i)}}, \sigma'_{V_{(i)}}\rangle \bigg] \label{eq:sig_and_orderstatistics_generalised}
\end{align}
\end{lemma}

\begin{proof}
\par It is known \cite{hambly2010uniqueness} that 
\begin{align} 
    n! S(\gamma)^n=\mathbb{E}\left[\bigotimes_{i=1}^n \gamma'_{U_{(i)}}\right] \label{eq:kth_order_sig_as_uniform}
\end{align}
\par Hence, for an orthonormal basis $\left\{e_\rho\right\}_{\rho \in \{1,..,n\}}$ of $\mathbb{R}^n$, we have
\begin{align}
    (n)^2\langle S(\gamma)^n, S(\sigma)^n \rangle &\overset{\eqref{eq:kth_order_sig_as_uniform}}{=} \mathbb{E} \left[ \prod_{i=1}^n \langle \gamma'_{U_{(i)}}, \sigma'_{V_{(i)}} \rangle\right]
\end{align}
For more details, we invite the reader to follow the arguments in \cite{hambly2010uniqueness}.
\end{proof}

\par This result states that the inner product between signatures of deterministic paths can be represented statistically through the mean of the product of $\langle \gamma'_{U_{(i)}}, \sigma'_{V_{(i)}}\rangle$. Observe that for unit-speed curves, the inner products $\langle \gamma'_{U_{(i)}}, \sigma'_{V_{(i)}}\rangle$ only encode the information on the angles $\Theta_i$ between the vectors $\gamma'_{U_{(i)}}$ and  $\sigma'_{V_{(i)}}$, i.e.
\begin{align}
    \langle \gamma'_{U_{(i)}}, \sigma'_{V_{(i)}}\rangle=\cos\left(\Theta_i\right):=\cos\left(\angle(\gamma'_{U_{(i)}},\sigma'_{V_{(i)}})\right).
\end{align}
In this case, also observe that the angles $\Theta_i$ can be exactly recovered from the norm of the difference between the two above random variables,
\begin{align}
    \Theta_i=\cos^{-1}\left(\langle \gamma'_{U_{(i)}}, \sigma'_{V_{(i)}}\rangle\right)=\cos^{-1}\left(1-\frac{1}{2}\left\vert\gamma'_{U_{(i)}} - \sigma'_{V_{(i)}} \right\vert^2\right). \label{def:angle_between_paths_formula}
\end{align}

\begin{proposition}[Inner product as a probabilistic expression]\label{prop:inner_product_sigs_as_sums}
Suppose that $\gamma,\sigma:\left[  0,1\right]  \rightarrow \mathbb{R}^d$ are two
absolutely continuous curves such that $|\gamma_{t}^{\prime}|=\left\vert
\sigma_{t}^{\prime}\right\vert =1$ for almost every \thinspace$t\in\left[
0,1\right]  .$ Then for every $n\in%
\mathbb{N}
$ we have
\begin{equation}\label{identity}
\begin{split}
\langle S(\gamma)^{n},S(\sigma)^{n}\rangle=&\frac{1}{(n!)^{2}}\mathbb{E}\left[
\exp\left(  -\sum_{k=1}^{\infty}\frac{1}{k2^k}\sum
_{i=1}^{n}\left\vert \gamma_{U_{\left(  i\right)  }}^{\prime}-\sigma
_{V_{\left(  i\right)  }}^{\prime}\right\vert ^{2k}\right)\id_{\{\max \Theta_i<\frac{\pi}{2}\}}  \right]\\
&+\frac{1}{(n!)^{2}}\mathbb{E}\left[
\id_{\{\max \Theta_i\geq\frac{\pi}{2}\}}\prod_{i=1}^{n} \langle \gamma'_{U_{(i)}}, \sigma'_{V_{(i)}}\rangle  \right],%
\end{split}
\end{equation}
where $\Theta_{i}$ in $\left[  0,\pi\right]  $ is defined by $\,\left\langle
\gamma_{U_{\left(  i\right)  }}^{\prime},\sigma_{V_{\left(  i\right)  }%
}^{\prime}\right\rangle \,=\cos\Theta_{i},$ for $i=1,\dots,n$, and $\id_A$ denotes the indicator function on a set $A$.
\end{proposition}
\begin{proof}
One the set $\left\{\max_{i=1,\dots,n}\Theta_i<\frac{\pi}{2}\right\}$, it holds that
\begin{equation*}
    \frac{1}{2}\left\vert\gamma'_{U_{(i)}} - \sigma'_{V_{(i)}} \right\vert^2<1,\quad\text{for }i=1,\dots,n.
\end{equation*}
Consequently, on this set, we have the expansion
\begin{align*}
    \log\left(\prod_{i=1}^{n} \langle \gamma'_{U_{(i)}}, \sigma'_{V_{(i)}}\rangle\right)&=\sum_{i=1}^n\log\left(1- \frac{1}{2}\left\vert\gamma'_{U_{(i)}} - \sigma'_{V_{(i)}} \right\vert^2\right)\\
    &=-\sum_{i=1}^n\sum_{k=1}^\infty \frac{1}{k}\left(\frac{1}{2}\left\vert\gamma'_{U_{(i)}} - \sigma'_{V_{(i)}} \right\vert^{2}\right)^{k}\\
    &=-\sum_{k=1}^\infty\frac{1}{k2^k}\sum_{i=1}^n\left\vert\gamma'_{U_{(i)}} - \sigma'_{V_{(i)}} \right\vert^{2k},
\end{align*}
where the second equality utilises the standard Talyor expansion for $\log(1-x)$, which is valid for $|x|<1$. Combining this expansion with \cref{lem:sig_and_order_stat_augmented} yields \cref{identity}.
\end{proof}
The following two subsections leverage the preceding probabilistic expression to attain upper and lower bounds for $\norm{S(\gamma)^n}$. We first collate the recurrent objects and assumptions that will be used in several subsequent arguments.

\begin{assumptions}[Standing assumptions]\label{def:usual_definitions}
Let $(\Omega, \mathcal{F}, \mathbb{P})$ be a probability space. The \textbf{standing assumptions} will refer to the following set of recurring assumptions and definitions,
\begin{itemize}
\item[(i)] Let $\gamma:[0,1]\rightarrow \mathbb{R}^d$ be a twice-continuously-differentiable unit-speed path such that the map $s \mapsto |\gamma_s''|$ is non-vanishing, and is differentiable with bounded derivative. Let $F:[a,b]\rightarrow [0,1]$ be its associated cumulative distribution as prescribed in \cref{lem:framing_lemma}.

\item[(ii)] Let $\{U_{(i)}\}_{i \in \{1,\dots,n\}}$ and $\{V_{(i)}\}_{i \in \{1,\dots,n\}}$ be two i.i.d. collections of uniform order statistics (\cref{def:uniform_order_statistics_samples}). Define the collections $X_{(i)}$ and $Y_{(i)}$ as
\begin{align}
    X_{(i)}:=F^{-1}(U_{(i)}), \quad Y_{(i)}:=F^{-1}(V_{(i)}), \quad i = 1,\dots,n.
\end{align}
\item[(iii)] Let $\mu_n^X$ and $\mu_n^Y$ be their respective empirical distributions defined by
\begin{align}
    \mu_n^X:=\frac{1}{n}\sum_{i=1}^n \delta_{X_{(i)}}, \quad \mu_n^Y:=\frac{1}{n}\sum_{i=1}^n \delta_{Y_{(i)}}.
\end{align}
\item[(iv)] Let $\mu$ and $f$ respectively be the probability measure and density distribution function associated with the cumulative distribution function $F$.
\end{itemize}
\end{assumptions}
To conclude this subsection we recall the definition of $p-$Wasserstein distances and a characterisation of distances between empirical measures.
\begin{definition}[Wasserstein Distance]
    Let $\mu$ and $\nu$ be probability measures supported on $\mathbb{R}$. The $p$\textsuperscript{th} Wasserstein distance, $W_p(\mu,\nu)$, between $\mu$ and $\nu$ is defined by
    \begin{equation}
    W_p^p(\mu,\nu)=\inf_{\pi\in\Gamma(\mu,\nu)}\int_\mathbb{\R}\int_{\mathbb{R}} \left\vert x-y\right\vert^p d\pi(x,y),
    \end{equation}
    where $\Gamma(\mu,\nu)$ denotes the set of all couplings of $\mu$ and $\nu$
\end{definition}
In the case that $\mu$ and $\nu$ are both discrete measures, then the infimum is explicit.
\begin{lemma}[Discrete characterisation of $p$-Wasserstein distance; Lemma 4.2 in \cite{bobkov2014one}] \label{lem:discrete_charac_wasserstein}
\par Let $\{A_j\}_{j\in\{1,..,n\}}$ and $\{B_j\}_{j\in\{1,..,n\}}$ be two samples of i.i.d. random variables. Denote by $\{A_{(j)}\}_{j\in\{1,..,n\}}$ and $\,\{B_{\left( j\right) }\}_{j\in\{1,..,n\}}$ their respective order statistics and let $\mu_{n}^{A}:=\frac{1}{n}\sum_{i=1}^{n}\delta _{A_{i}}$
and $\mu _{n}^{B}=\frac{1}{n}\sum_{i=1}^{n}\delta _{B_{i}}$ be their associated empirical
probability measures. Then we have%
\begin{equation*}
\sum_{j=1}^{n}\left( A_{\left( j\right) }-B_{\left( j\right) }\right)
^{p}=nW_{p}^{p}\left( \mu _{n}^{A},\mu _{n}^{B}\right) ,
\end{equation*}
where $W_{p}$ is the $p$-Wasserstein distance.
\end{lemma}
\subsection{Lower bound on $\norm{S(\gamma)^n}$ in terms of Wasserstein distances}\label{subsec:lower_b_sig_in_terms_of_wasserstein}
We now use the probabilistic expression of the signature inner product (\cref{prop:inner_product_sigs_as_sums}) when $\sigma=\gamma$, and derives a lower-bound on $\norm{S(\gamma)}$ in terms of Wasserstein distances.

\begin{proposition}[Lower bound on $\norm{S(\gamma)}$ in terms of Wasserstein distances]\label{prop:lower_bound_sig_wasserstein}
Let $\gamma,\; \mu_n^X,\;\mu_n^Y$ be as defined in \cref{def:usual_definitions}. Then
\begin{equation}\label{lower bound}
\begin{split}
(n!)^2\norm{ S\left(  \gamma\right)  ^{n}}^2
\geq&\mathbb{E}\left[  \exp\left(  -n\sum_{k=1}^{\infty}\frac{1}{k2^k}W_{2k}^{2k}(\mu_{n}^{X},\mu_{n}^{Y})\right)\id_{\left\{\max_{i=1,...,n} \left\vert X_{(i)}-Y_{(i)}\right\vert<1\right\}}
\right]\\
&-\mathbb{P}\left(  \max_{i=1,...,n}\left\vert X_{\left(  i\right)  }-Y_{\left(
i\right)  }\right\vert \geq 1\right),
\end{split}
\end{equation}
for every $n\in \mathbb{N}$.
\end{proposition}

\begin{proof}
We use the fact that
\begin{equation}
|\gamma_{U_{(i)}}^{\prime}-\gamma_{V_{(i)}}^{\prime}|=|g(X_{\left(  i\right)
})-g(Y_{\left(  i\right)  })| \label{inverse identity}%
\end{equation}
where $g:=\gamma^{\prime}\circ F$. The assumptions on $\gamma$ give that $g$ is once continuously differentiable and so the mean value inequality may be employed to see that
\begin{equation}\label{eq:norm_1}
    |\gamma_{U_{(i)}}^{\prime}-\gamma_{V_{(i)}}^{\prime}|\leq|X_{\left(  i\right)
}-Y_{\left(  i\right)  }|\cdot|g^{\prime}(\xi_{i})|\text{ for some }\xi_{i}%
\in\lbrack X_{\left(  i\right)  },Y_{\left(  i\right)  }]\text{ for
}i=1,\dots,n.
\end{equation}
Furthermore, as $F$ satisfies the integral equation of \cref{lem:framing_lemma} we
have that
\[
|g^{\prime}(s)|=|\gamma_{F(s)}^{\prime\prime}|\cdot|F^{\prime}(s)|=1,\;\forall
s\in[a,b]
\]

By applying \cref{lem:discrete_charac_wasserstein} we learn that
\begin{equation}
\sum_{i=1}^{n}\left\vert \gamma_{U_{\left(  i\right)  }}^{\prime}%
-\gamma_{V_{\left(  i\right)  }}^{\prime}\right\vert ^{2k}\leq\sum_{i=1}%
^{n}\left\vert X_{\left(  i\right)  }-Y_{\left(  i\right)  }\right\vert
^{2k}=nW_{2k}^{2k}(\mu_{n}^{X},\mu_{n}^{Y}). \label{W est}%
\end{equation}
Observe that $\max_{i=1,\dots,n} \left\vert X_{(i)}-Y_{(i)}\right\vert<1$ is a strictly stronger condition than $\max_{i=1,\dots,n} \Theta_i<\frac{\pi}{2}$, and that the product in the second term on the right-hand-side of \cref{identity} may be lower bounded by $-1$ by the Cauchy-Schwarz inequality. The lower bound \cref{lower bound} then follows by combining this observation, \cref{identity}, and
\cref{W est}.
\end{proof}
\begin{remark}\label{rem:ind_inside_exp}
We can also multiply the sum inside the exponential term in \cref{lower bound} by the indicator function on the set $\{\max_{i=1,...,n} \left\vert X_{(i)}-Y_{(i)}\right\vert<1\}$ without changing the random variable inside the expectation. Doing so will prove convenient in the proof of our main result \cref{thm:generalisation_hambly_lyons}.
\end{remark}
\subsection{Upper bound on $\norm{S(\gamma)^n}$ in terms of Wasserstein distances}\label{subsec:upper_b_sig_in_terms_of_wasserstein}
Similarly to \cref{subsec:lower_b_sig_in_terms_of_wasserstein}, we use \cref{prop:inner_product_sigs_as_sums} to derive an upper bound on for $\norm{S(\gamma)^n}$ in terms of a series of Wasserstein distances.
\begin{proposition}[Upper bound on $\norm{S(\gamma)}$ in terms of Wasserstein distances]\label{prop:upper_bound}
    Let $\gamma,\; \mu_n^X,\;\mu_n^Y$ be as defined in \cref{def:usual_definitions}. Then there exists an $0<\varepsilon^\prime\leq 1$ so that for every $n\in \mathbb{N}$ and $0<\varepsilon<\varepsilon^\prime$
    \begin{equation}\label{eq:upper_bound}
    \begin{split}
    (n!)^2\norm{ S\left(  \gamma\right)  ^{n}}^2
    \leq&\mathbb{E}\left[  \exp\left(  -n\sum_{k=1}^{\infty}\frac{(1-\phi(\varepsilon))^k}{k2^k}W_{2k}^{2k}(\mu_{n}^{X},\mu_{n}^{Y})\right)\id_{\left\{\max_{i=1,...,n} \left\vert X_{(i)}-Y_{(i)}\right\vert<\varepsilon\right\}}
    \right]\\
    &+\mathbb{P}\left(  \max_{i=1,...,n}\left\vert X_{\left(  i\right)  }-Y_{\left(
    i\right)  }\right\vert \geq \varepsilon\right),
    \end{split}
    \end{equation}
    where $\phi\left(  u\right)  :=\sup_{\left\vert s-t\right\vert <u}\left\vert
g^{\prime}\left(  t\right)  -g^{\prime}\left(  s\right)  \right\vert $ is the modulus of continuity of the derivative of $g:=\gamma^\prime\circ F$.
\end{proposition}
\begin{proof}
    An application of the fundamental theorem of calculus and \cref{eq:norm_1} gives
    \begin{align*}
        \left\vert\gamma_{U_{(i)}}^{\prime}-\gamma_{V_{(i)}}^{\prime}\right\vert^{2}&=\left\vert g(X_{\left(  i\right)
        })-g(Y_{\left(  i\right)  })\right\vert^2\\
        &=\left\langle \int_{X_{\left(  i\right)}}^{Y_{\left(  i\right)}}g^{\prime}(r)dr,\int_{X_{\left(  i\right)}}^{Y_{\left(  i\right)}}g^\prime(s)ds\right\rangle\\
        &=\int_{\left[X_{\left(  i\right)},Y_{\left(  i\right)}\right]^2}\left\langle g^\prime(r),g^\prime(s)\right\rangle -\left\vert g^\prime(r)\right\vert^2 + 1dr ds\\
        &= \left\vert X_{\left(  i\right)} - Y_{\left(  i\right)}\right\vert^2 + \int_{\left[X_{\left(  i\right)},Y_{\left(  i\right)}\right]^2}\left\langle g^\prime(r),g^\prime(s)-g^\prime(r)\right\rangle drds.
    \end{align*}
    And hence, by applying the Cauchy-Schwarz inequality to the integrand in the last line
    \begin{equation}
        \left\vert \left\vert\gamma_{U_{(i)}}^{\prime}-\gamma_{V_{(i)}}^{\prime}\right\vert^{2}-\left\vert
        X_{\left(  i\right)  }-Y_{\left(  i\right)  }\right\vert  ^{2}\right\vert
        \leq\left\vert  X_{\left(  i\right)  }-Y_{\left(  i\right)  }\right\vert  ^{2}%
        \phi\left(  \left\vert X_{\left(  i\right)  }-Y_{\left(  i\right)
        }\right\vert \right),
    \end{equation}
    and so, in particular
    \begin{equation}\label{eq:term_wise_ub}
        \left\vert\gamma_{U_{(i)}}^{\prime}-\gamma_{V_{(i)}}^{\prime}\right\vert^{2}\geq \left(1-\phi\left(  \left\vert X_{\left(  i\right)  }-Y_{\left(  i\right)
        }\right\vert \right)\right)\left\vert  X_{\left(  i\right)  }-Y_{\left(  i\right)  }\right\vert  ^{2}.
    \end{equation}
    To ensure that the series inside the exponential in \cref{eq:upper_bound} is finite, we note that $\phi$ is a modulus of continuity, and so there exists some $0<\varepsilon^\prime\leq 1$ for which $\phi(\varepsilon)<1$ for any $\varepsilon<\varepsilon^\prime$. Now, the product in the second term on the right-hand-side of \cref{identity} may be upper bounded by $1$ by the Cauchy-Schwarz inequality, so it follows from \cref{eq:term_wise_ub}, \cref{lem:discrete_charac_wasserstein} and \cref{prop:inner_product_sigs_as_sums} that \cref{eq:upper_bound} holds provided $\varepsilon<\varepsilon^\prime$.
\end{proof}
\subsection{Generalising the Hambly-Lyons Limit Theorem}\label{subsec:generalising_hambly_lyons}
Combined, \cref{prop:lower_bound_sig_wasserstein,prop:upper_bound} provide lower and upper bounds for $\norm{S(\gamma)^n}$ in terms of a series of Wasserstein distances. What remains is to show that the lower bound converges to the square of the Hambly-Lyons limit as $n\to\infty$, and that the same applies to the upper bound when taking $n\to\infty$ and then $\varepsilon\to 0$. The following pair of lemmas provide the necessary results for this conclusion.
\begin{lemma}\label{lem:prob_to_zero}
    Let $\mu$, $\mu_n^X$, and $\mu_n^Y$ be as in the standing assumptions, then for any $\varepsilon>0$
    \begin{equation}
        \lim_{n\to\infty} \mathbb{P}\left(  \max_{i=1,...,n}\left\vert X_{\left(  i\right)  }-Y_{\left(
    i\right)  }\right\vert \geq \varepsilon\right)=0.
    \end{equation}
\end{lemma}
\begin{proof}
Using mean value and inverse function Theorems, we may deduce that
\begin{equation}\label{eq:term_wise_equal}
X_{\left(  i\right)  }-Y_{\left(  i\right)  }=F^{-1}\left(  U_{\left(
i\right)  }\right)  -F^{-1}\left(  V_{\left(  i\right)  }\right)  =\left(
F^{-1}\right)  ^{\prime}\left(  \eta_{i}\right)  \left(  U_{\left(  i\right)
}-V_{\left(  i\right)  }\right)  =\left\vert\gamma_{\eta_{i}%
}^{\prime\prime}\right\vert \left(  U_{\left(  i\right)
}-V_{\left(  i\right)  }\right)
\end{equation}
for some $\eta_{i}\in\left[  U_{\left(  i\right)  },V_{\left(  i\right)
}\right]$. An application of Markov's inequality gives that
\begin{equation*}
    \mathbb{P}\left(  \max_{i=1,...,n}\left\vert X_{\left(  i\right)  }-Y_{\left(
i\right)  }\right\vert \geq\varepsilon\right)\leq \varepsilon^{-1}\norm{\gamma^{\prime\prime}}_{\infty}\mathbb{E}\left(  \max_{i=1,...,n}\left\vert U_{\left(  i\right)  }-V_{\left(
i\right)  }\right\vert\right)\leq 2C\varepsilon^{-1}\norm{\gamma^{\prime\prime}}_{\infty}\frac{1}{\sqrt{n}},
\end{equation*}
for some absolute constant $C>0$. The first inequality utilises \cref{eq:term_wise_equal}, and the second is due to Theorem 4.9 of \cite{bobkov2014one}. Taking the limit as $n\to\infty$ in the above inequality concludes the proof.
\end{proof}
\begin{lemma}\label{lem:L_1_convergence}
Let $\mu$, $\mu_n^X$, and $\mu_n^Y$ be as in the standing assumptions and $0<a,\varepsilon\leq 1$, then
\begin{equation}
    \Delta_n := n\id_{\left\{\max_{i=1,...,n} \left\vert X_{(i)}-Y_{(i)}\right\vert<\varepsilon\right\}}\sum_{k=2}^{\infty}\frac{a^k}{k2^k}W_{2k}^{2k}(\mu_{n}^{X},\mu_{n}^{Y})\rightarrow 0
\end{equation}
in $L_1$ as $n\to\infty$.
\end{lemma}
\begin{proof}
    We have that
    \begin{align*}
        \mathbb{E}\left[\Delta_n\right]&=n\mathbb{E}\left[\id_{\left\{\max_{i=1,...,n} \left\vert X_{(i)}-Y_{(i)}\right\vert<\varepsilon\right\}}\sum_{k=2}^{\infty}\frac{a^k}{k2^k}W_{2k}^{2k}(\mu_{n}^{X},\mu_{n}^{Y})\right]\\
        &=\mathbb{E}\left[\id_{\left\{\max_{i=1,...,n} \left\vert X_{(i)}-Y_{(i)}\right\vert<\varepsilon\right\}}\sum_{k=2}^{\infty}\frac{a^k}{k2^k}\sum_{i=1}
        ^{n}\left\vert X_{\left(  i\right)  }-Y_{\left(  i\right)  }\right\vert
        ^{2k}\right]\\
        &\leq \mathbb{E}\left[\sum_{k=2}^{\infty}\frac{a^k}{k2^k}\sum_{i=1}
        ^{n}\left\vert X_{\left(  i\right)  }-Y_{\left(  i\right)  }\right\vert
        ^{4}\right]\\
        &\leq \mathbb{E}\left[nW_4^4(\mu_n^X,\mu_n^Y)\sum_{k=1}^{\infty}\frac{a^k}{k2^k}\right]\\
        &=\mathbb{E}\left[nW_4^4(\mu_n^X,\mu_n^Y)\right]\log\left(\frac{2}{2-a}\right).
    \end{align*}
    Using the following inequality (Section 4.1, \cite{bobkov2014one}) for $p\in\mathbb{N}\setminus\{0\}$,
    \begin{equation*}
        \mathbb{E}\left[W_p^p(\mu_n^X,\mu_n^Y)\right]\leq 2^p\mathbb{E}\left[W_p^p(\mu_n^X,\mu)\right],
    \end{equation*}
    together with the bounds from Theorem 5.3 in \cite{bobkov2014one} gives
    \begin{equation*}
        \mathbb{E}\left[W_4^4(\mu_n^X,\mu_n^Y)\right]\leq 2^4 \left(\frac{20}{\sqrt{n+2}}\right)^4J_4(\mu),
    \end{equation*}
    where $J_4(\mu)$ is finite since $\mu$ is compactly supported with density $f$ bounded below by \cref{lem:existence_sol_framing}. It follows that
    \begin{equation*}
        \lim_{n\to\infty}\mathbb{E}\left[\Delta_n\right]\leq \lim_{n\to\infty}40^4\log\left(\frac{2}{2-a}\right)\frac{n}{(n+2)^2}J_4(\mu) =0.
    \end{equation*}
\end{proof}
We are now ready to prove our main result.
\begin{theorem}[Generalisation of the Hambly-Lyons Limit Theorem]\label{thm:generalisation_hambly_lyons}
\par Let $\gamma$, $\mu$, $\mu_n^X$, and $\mu_n^Y$ be as in the standing \cref{def:usual_definitions}. Then,
\begin{equation}\label{eq:main_result}
    \lim_{n\rightarrow \infty} n!\norm{S(\gamma)^n}=\mathbb{E}\left[\exp\left(-\int_0^1 (B_{0,s})^2|\gamma_s''|^2ds\right)\right]^{1/2}.
\end{equation}
\end{theorem}
\begin{proof}
   \par For $0<\varepsilon\leq 1$, define $A_n^\varepsilon$ to be the set $\{\max_{i=1,\dots,n}|X_{(i)}-Y_{(i)}|<\varepsilon\}$. Then, by \cref{lem:prob_to_zero}, $\id_{A_n^\varepsilon}$ converges to $1$ in probability, and that  $\mathbb{P}\left((A_n^\varepsilon)^c\right)$ converges to $0$. It follows from \cref{cor:BGU_empirical,cor:coinicde_hamlby_lyons_bgu}, and Slutsky's Theorem that 
    \begin{equation}\label{eq:indicator_conv_bgu}
        \frac{n}{2}\id_{A_n^\varepsilon}W_2^2(\mu_n^X,\mu_n^Y)\to \int_0^1 (B_{0,s})^2|\gamma_s''|^2ds
    \end{equation}
    in distribution as $n\to\infty$. As such, by \cref{eq:indicator_conv_bgu}, \cref{lem:L_1_convergence}, and another application of Slutsky's Theorem we obtain the following convergence in distribution
    \begin{equation}\label{eq:conv_1}
        n\id_{A_n^\varepsilon}\sum_{k=1}^{\infty}\frac{1}{k2^k}W_{2k}^{2k}(\mu_{n}^{X},\mu_{n}^{Y})\to \int_0^1 (B_{0,s})^2|\gamma_s''|^2ds.
    \end{equation}
    For the lower bound, the combining of \cref{eq:conv_1} with \cref{prop:lower_bound_sig_wasserstein}, \cref{rem:ind_inside_exp}, \cref{lem:prob_to_zero}, yet another application of Slutsky's Theorem, and the Continuous Mapping Theorem, results in the lower bound
    \begin{equation}\label{eq:final_lower_bound}
        \lim_{n\to\infty} (n!)^2\norm{S(\gamma)^n}^2\geq \mathbb{E}\left[\exp\left(-\int_0^1 (B_{0,s})^2|\gamma_s''|^2ds\right)\right].
    \end{equation}
    Similar analysis for the upper bound from \cref{prop:upper_bound} shows that
    \begin{equation}\label{eq:final_upper_bound}
        \lim_{n\to\infty} (n!)^2\norm{S(\gamma)^n}^2\leq \mathbb{E}\left[\exp\left(-(1-\phi(\varepsilon))^2\int_0^1 (B_{0,s})^2|\gamma_s''|^2ds\right)\right],
    \end{equation}
    for suitably small $\varepsilon$. Using the fact that $\phi$ is continuous at zero, we may take the limit in the preceding as $\varepsilon\to 0$ and combine it with the lower bound \cref{eq:final_lower_bound} to conclude that \cref{eq:main_result} holds.
\end{proof}

\par This generalisation allows us to compute the limit of the signature norm for curves that are not $C^3$.

\begin{example}[Integral of semi-circular curves]
\par Let $R>\frac{1}{2\pi}$ be a constant. Let $\gamma_+,\gamma_-:[0,\pi R]\rightarrow \mathbb{R}^2$ denote two planar semi-circular curves in the plane of radius $R$, rotating clockwise and anti-clockwise respectively, and defined as
\begin{align}
    \gamma_\pm(s):=R\left(\cos\left[\mp\; \frac{s}{R}\right],\sin\left[\mp\; \frac{s}{R}\right]\right) .
\end{align}
\par Consider a concatenation $\psi$ of these two curves, say $\psi:=\gamma_+ * \gamma_-:[0,2\pi R]\rightarrow \mathbb{R}^2$. It is easy to see that $\psi$ is differentiable everywhere and that its second-derivative exists everywhere except at point $t=\frac{\pi}{R}$. Define now the integral $\Psi:[0,2\pi R]\rightarrow \mathbb{R}^2$ as
\begin{align}
    \Psi(t):= \frac{1}{R}\int_{0}^t \psi(s)ds.
\end{align}
Then, $\Psi|_{[0,1]}$ satisfies the condition for \cref{thm:generalisation_hambly_lyons} but is not $C^3$. Indeed,
\begin{itemize}
    \item[(i)] Because $\psi$ is continuously-differentiable, $\Psi$ is twice-continuously-differentiable by the fundamental theorem of calculus. 
    \item[(ii)] As $\Psi''(s)=\psi'(s)$, the map $s\mapsto |\Psi''(s)|=\frac{1}{R}$ is non-vanishing, and is differentiable with bounded derivative.
    \item[(iii)] Because $\psi$ is not twice-differentiable everywhere, $\Psi$ does not admit a third-derivative everywhere and is therefore not $C^3$.
\end{itemize}
\end{example}

\begin{figure}[H]
\begin{center}
  \includegraphics[width=12cm]{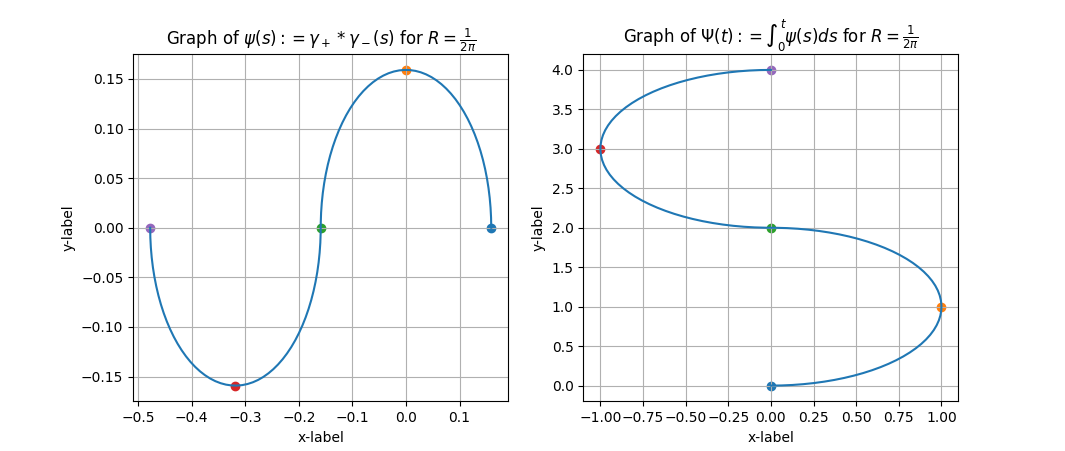}
  \caption{Graph of $\phi$ and $\Phi$ for $R=\frac{1}{2\pi}$. The blue, orange, green, red, and purple dots denote evaluation at time $t \in \{0,\frac{1}{4}, \frac{1}{2}, \frac{3}{4}, 1\}$ respectively.}
  \end{center}
\end{figure}


\subsection{Computing the Hambly-Lyons limit explicitly}\label{subsection:computing_hambly_lyons}

\par Finally, we propose a way to practically compute the limit presented in \cref{thm:generalisation_hambly_lyons}. Using the work of Yor and Revuz on Bessel bridges \cite{revuz2013continuous}, we relate the computation of the expected integral $c(\gamma)$ to the solving of a second order distributional differential equation.

\begin{lemma}[High order term and curvature]\label{lem:yor_revuz_integral_into_dde}
Let $V$ be a finite dimensional inner product space. Suppose that
$\gamma:\left[  0,l\right]  \rightarrow V$ is parameterised at unit speed and
is three-times continuously differentiable (under this parameterisation). Let
$\mu$ denote the finite Borel measure on $%
\mathbb{R}
_{+}$ which is absolutely continuous with respect to the Lebesgue measure
$\lambda=\lambda_{%
\mathbb{R}
_{+}}$ with density given by $\frac{d\mu}{d\lambda}\left(  t\right)
=l^{2}\left\vert \gamma_{tl}^{\prime\prime}\right\vert ^{2}1_{[0,1]}\left(
t\right)  $ a.e. Then,
\begin{itemize}
\item[(i)] There exists a unique continuous function $\phi:\left[
0,1\right]  \rightarrow\mathbb{R}$ which solves, in the distributional sense,
the second-order differential equation
\begin{equation}
\psi^{\prime\prime}=2\psi\mu,\text{with }\psi_{0}=0,\psi_{0}^{\prime
}=1.\label{ode}%
\end{equation}

\item[(ii)] This solution satisfies \[
\psi_1^{-l/4}  =\mathbb{E}\left[  \exp\left(  -\int_{0}^{1}\left(
B_{0,t}\right)  ^{2}\mu\left(  dt\right)  \right)  \right]  ^{1/2}=c(\gamma).
\]
\end{itemize}
\end{lemma}

\begin{proof}
Define the unit-speed curve $\tilde{\gamma}_{t}:=\frac{1}{l}\gamma_{lt}$ over
$\left[  0,1\right]  $ and let the terms in its signature be given by $\left(
1,S(\tilde{\gamma})^{1},..,S(\tilde{\gamma})^{n},..\right).$ It suffices to prove the result
under the assumption that $l=1$ since the general result can then be recovered
by noticing
\[
\left(  n!\norm{ S(\gamma)^{n}} \right)^{1/n}=l\left(  n!\norm{S(\tilde{\gamma})^n
}\right)^{1/n}\sim l+\frac{l}{n}\log c\left(\tfrac{1}{l}\gamma\right)  +o\left(  \frac{1}{n}\right)  ,
\]
so that $c\left(  \gamma\right)  =c\left(\tfrac{1}{l}\gamma\right)  ^{l}.$ We therefore assume
that $l=1,$ and seek to prove that $c\left(  \gamma\right)  =\psi_{1}^{-1/4}.$ To
this end, we first observe that the unit-speed parameterisation of $\gamma$
gives that $\left\langle \gamma_{t}^{\prime},\gamma_{t}^{\prime\prime\prime
}\right\rangle =-\left\vert \gamma_{t}^{\prime\prime}\right\vert ^{2}$ for
every \thinspace$t$ in $\left[  0,1\right]  .$ When used together with
\cref{HL id} this gives that%
\[
c\left(  1\right)  =\mathbb{E}\left[  \exp\left(  -\int_{0}^{1}\left(
B_{0,t}\right)  ^{2}\mu\left(  dt\right)  \right)  \right]  ^{1/2}.
\]
By noticing that $\left(  B_{\cdot}^{0}\right)  ^{2}$ is a squared Bessel
bridge of dimension $1$ which returns to zero at time $1$, we can then use
Theorem 3.2 of \cite{revuz2013continuous} to see that
\begin{equation}
c\left(  \gamma\right)  =\left(  \phi_{1}\int_{0}^{1}\frac{1}{\phi_{t}^{2}%
}dt\right)  ^{-1/4},\label{c1}%
\end{equation}
where $\phi$ is the unique continuous positive function solving the
distributional equation $\phi^{\prime\prime}=2\phi\mu,$ with $\phi_{0}%
=1,\;\phi_{0}^{\prime}=0.$. Using
Exercise 1.34 in \cite{revuz2013continuous}, we know that the function $\psi$ in the
statement and $\phi$ are related by
\[
\psi_{t}=\phi_{t}\int_{0}^{t}\frac{1}{\phi_{s}^{2}}ds,
\]
and therefore \cref{c1} becomes $c\left(  \gamma\right)  =\psi_{1}^{-1/4}.$
\end{proof}

\par This result reduces the computation of $c(\gamma)$ for a planar circle to the solving of a simple differential equation.

\begin{example}[The planar circle]\label{examp:planar_circle}
\par Let $\gamma_{t}=\frac{1}{2\pi}\left(  \sin(2\pi
t),\cos(2\pi t)\right)  $ for $t$ in $\left[  0,1\right]  $, then $\gamma$ is a
smooth unit-speed curve of length $l=1$ and $\left\vert \gamma_{t}%
^{\prime\prime}\right\vert =2\pi$ so that $\mu\left(  dt\right)
=2\pi1_{[0,1]}\left(  t\right)  dt$. By solving the differential equation
\cref{ode} we find that
\begin{align}
    c\left(  \gamma\right)  =\left(  \frac{2\sqrt{\pi}}{\sinh2\sqrt{\pi}}\right)
    ^{1/4}.
\end{align}
\end{example}

\begin{remark}[Role of $c(\gamma)$ in the asymptotics of the norm]
\par By starting from \cref{weak} and
\cref{HL id}, we notice the latter statement can be rewritten as a
statement of asymptotic equivalence, namely
\begin{equation}
\left(  n!\norm{ S(\gamma)^{n}} _{\text{HS}%
}\right)  ^{1/n}\sim l+\frac{1}{n}\log c\left(  \gamma\right)  +o\left(  \frac
{1}{n}\right)  \text{ as }n\rightarrow\infty\text{ ,} \label{eq:original_asymptotics_hilbert_schmidt} \end{equation}
where we write $\norm{\cdot} _{\text{HS
}}$ for the moment to emphasise the dependence of this expansion on the choice
of tensor norm. By contrast for the projective tensor norm it follows from
\cref{strong} that we have
\[
\left(  n!\norm{ S(\gamma)^{n}} _{\text{proj}%
}\right)  ^{1/n}\sim l+o\left(  \frac{1}{n}\right)  \text{ as }n\rightarrow
\infty.
\]
\par When written in this way, the $c(\gamma)$ limit has the
interpretation of being the second term in the asymptotic expansion of
$(n!\norm{ S(\gamma)^{n}} _{\text{HS}})^{1/n}$ as
$n\rightarrow\infty.$ Natural questions would be to explore the higher order terms in these asymptotic expansions, and to relate them to geometric features of the underlying curve.
\end{remark}

\section*{Funding}
Thomas Cass has been supported by the EPSRC Programme Grant EP/S026347/1. Remy Messadene and William F. Turner have been supported by the EPSRC Centre for Doctoral Training in Mathematics of Random Systems: Analysis, Modelling and Simulation (EP/S023925/1).

\printbibliography

\end{document}